\documentclass[12pt]{article}
\usepackage{mathrsfs}
\usepackage{stmaryrd}
\usepackage{amsfonts,amsmath,amssymb,amscd}
\usepackage{shadow}
\usepackage{hyperref}

\allowdisplaybreaks

\newtheorem{thm}{Theorem}[section]
\newtheorem{lem}[thm]{Lemma}

\newtheorem{conj}{Conjecture}[section]

\def\qed{\hfill \rule{4pt}{7pt}}

\def\pf{\noindent {\it{Proof.} \hskip 2pt}}

\numberwithin{equation}{section}

\begin{document}
\begin{center}
{\large\bf The  Log-Behavior of $\sqrt[n]{p(n)}$ and $\sqrt[n]{p(n)/n}$}
\end{center}

\begin{center}
William Y.C. Chen$^{1}$ and Ken Y. Zheng$^{2}$\\[8pt]
Center for Combinatorics, LPMC-TJKLC\\
Nankai Univercity\\
Tianjin 300071, P. R. China\\[6pt]
Email:$^{1}${\tt chenyc@tju.edu.cn}, $^{2}${\tt kenzheng@aliyun.com}
\end{center}
\begin{abstract}

Let $p(n)$ denote the partition function. Desalvo and Pak proved the log-concavity of $p(n)$ for $n>25$ and the
inequality $\frac{p(n-1)}{p(n)}\left(1+\frac{1}{n}\right)>\frac{p(n)}{p(n+1)}$
 for $n>1$. Let $r(n)=\sqrt[n]{p(n)/n}$ and $\Delta$ be the difference operator respect to $n$.
  Desalvo and Pak pointed out that their approach to
  proving the log-concavity of $p(n)$ may be employed to prove
   a conjecture of Sun on the log-convexity of
  $\{r(n)\}_{n\geq 61}$, as long as one finds an appropriate estimate of $\Delta^2 \log r(n-1)$. In this paper, we obtain  a lower bound for $\Delta^2\log r(n-1)$,  leading to
  a proof of this conjecture. From the log-convexity of $\{r(n)\}_{n\geq61}$ and $\{\sqrt[n]{n}\}_{n\geq4}$, we
  are led to a proof of another  conjecture of Sun on the log-convexity of $\{\sqrt[n]{p(n)}\}_{n\geq27}$.
 Furthermore, we show that $\lim\limits_{n \rightarrow +\infty}n^{\frac{5}{2}}\Delta^2\log\sqrt[n]{p(n)}=3\pi/\sqrt{24}$.
 Finally, by finding an upper bound of $\Delta^2 \log\sqrt[n-1]{p(n-1)}$, we prove an inequality on the ratio  $\frac{\sqrt[n-1]{p(n-1)}}{\sqrt[n]{p(n)}}$ analogous to the above  inequality on the ratio $\frac{p(n-1)}{p(n)}$.
\end{abstract}



\section{Introduction}

In this paper, we study the log-behavior of the
sequences $\sqrt[n]{p(n)}$ and $\sqrt[n]{p(n)/n}$, where $p(n)$ is the
number of partitions of $n$. Recently, by using the Hardy-Ramanujan-Rademacher formula of $p(n)$ (see \cite{Hardy1,Hardy,Rade}) and Lehmer's error bound (see \cite{Leh1,Leh2}),
Desalvo and Pak \cite{Pak} gave an estimate for $-\Delta^2 \log p(n-1)$, and then found an upper and lower bound for this estimate, finally proved that $p(n)$ is log-concave for $n>25$. They also proved the following inequality conjectured by Chen \cite{Chen}.

\begin{thm}
For $n>1$,
\begin{equation}
\frac{p(n-1)}{p(n)}\left(1+\frac{1}{n}\right)>\frac{p(n)}{p(n+1)}.
\end{equation}
\end{thm}

Desalvo and Pak \cite{Pak} showed  that
\begin{equation}\label{limdeltapn}
\lim\limits_{n \rightarrow
+\infty}-n^{\frac{3}{2}}\Delta^2\log p(n-1)
=\pi/\sqrt{24},
\end{equation}
and proposed the following conjecture.

\begin{conj}\label{dpc}
For $n\geq45$,
\begin{equation}\label{ineqdpc}
\frac{p(n-1)}{p(n)}\left(1+\frac{\pi}{\sqrt{24}n^{3/2}}\right)>\frac{p(n)}{p(n+1)}.
\end{equation}
\end{conj}

In view of \eqref{limdeltapn}, the coefficient $\frac{\pi}{\sqrt{24}}$
 in   \eqref{ineqdpc} is the best possible.
 Chen, Wang and Xie \cite{Chen1} proved the above conjecture
 by showing that  for $n\geq5000$,
 \begin{align}\label{in1}
 -\Delta^2 \log p(n-1)<\frac{24\pi}{(24n)^{3/2}}-
 \left(\frac{24\pi}{(24n)^{3/2}}\right)^2.
 \end{align}
 The proof of \eqref{in1} requires Desalvo and Pak's upper bound of $-\Delta^2 \log p(n-1)$ for $n\geq50$,
 \begin{align}\label{in2}
 \notag -\Delta^2 \log p(n-1)<&\frac{24\pi}{(24(n-1)-1)^{3/2}}
 +\frac{288\pi(-3+\pi\sqrt{24(n-1)-1})}{(24(n-1)-1)^{3/2}(-6+\pi\sqrt{24(n-1)-1})^2}\\[6pt]
&\qquad -\frac{864}{(24(n+1)-1)^2}
 +2e^{-\frac{\pi}{10}\sqrt{\frac{2n}{3}}}.
 \end{align}
 For $n\geq 50$, the upper bound in \eqref{in2} can be relaxed to
 \begin{align*}
 \frac{24\pi}{(24n)^{3/2}}-
 \left(\frac{24\pi}{(24n)^{3/2}}\right)^2
 -\frac{1}{n^2}+\frac{3}{n^{5/2}}+
 2e^{-\frac{\pi}{10}\sqrt{\frac{2n}{3}}}.
 \end{align*}
 By using the Lambert $\emph{W}$ function, it can be shown that
 $-\frac{1}{n^2}+\frac{3}{n^{5/2}}+
 2e^{-\frac{\pi}{10}\sqrt{\frac{2n}{3}}}<0$   when $n\geq5000$,
and therefore we arrive at the upper bound in the form of \eqref{in1}.
Let $r(n)=\sqrt[n]{p(n)/n}$. Desalvo and Pak also considered the log-behavior of $r(n)$. A positive sequence $\{a_n\}$ is log-convex if it satisfies that for $n\geq1$,
\[a_n^2-a_{n-1}a_{n+1}\leq0.\]  Conversely, a positive sequence $\{a_n\}$ is log-concave if it satisfies that for $n\geq1$,
\[a_n^2-a_{n-1}a_{n+1}\geq0. \]
Desalvo and Pak noticed that the log-convexity of $\{r(n)\}_{n\geq61}$ conjectured by Sun \cite{Sun}
 can be derived from an estimate for $\Delta^2 \log r(n-1)$,
 see \cite[Final Remark 7.7]{Pak}. They also remarked that
their approach to bounding
 $-\Delta^2 \log p(n-1)$ does not seem to
 apply to   $\Delta^2 \log r(n-1)$.
 In this paper, we obtain a lower bound for $\Delta^2 \log r(n-1)$, leading to a proof of the log-convexity of $\{r(n)\}_{n\geq61}$.

\begin{thm}\label{th1}
The sequence $\{r(n)\}_{n\geq61}$ is log-convex.
\end{thm}

The log-convexity of $\{r(n)\}_{n\geq61}$ implies the log-convexity of $\{\sqrt[n]{p(n)}\}_{n\geq27}$, because the
sequence $\{\sqrt[n]{n}\}_{n\geq4}$ is log-convex \cite{Sun}. It is known that $\lim \limits_{n \rightarrow +\infty}\sqrt[n]{p(n)}=1$. For a combinatorial
 proof of this fact, see Andrews \cite{Andrews1}. The log-convexity of $\{\sqrt[n]{p(n)}\}_{n\geq27}$ was conjectured
by Sun \cite{Sun}. He also proposed
the conjecture that $\{\sqrt[n]{p(n)}\}_{n\geq6}$ is strictly decreasing, which has been proved by Wang and Zhu \cite{Wang}. It is easy to see that the log-convexity of $\{\sqrt[n]{p(n)}\}_{n\geq27}$ implies the decreasing property.

 It should be noted that there is another approach to proving the log-convexity of $\{\sqrt[n]{p(n)}\}_{n\geq27}$. Chen, Guo and Wang \cite{Chen-a} introduced the notion of a ratio log-convex sequence and showed that  ratio
 log-convexity implies   log-convexity under
 an initial condition. A sequence $\{a_n\}_{n\geq k}$ is called ratio log-convex if $\{a_{n+1}/a_n\}_{n\geq k}$ is log-convex, or, equivalently, for $n\geq k$,
 \[\log a_{n+2}-3\log a_{n+1}+3\log a_{n}-\log a_{n-1}\geq0.\]
 Chen, Wang and Xie \cite{Chen1}
  showed that  that for any $r\geq1$, one can determine
   a number  $n(r)$ such that for $n>n(r)$,  $(-1)^{r-1}\Delta^r\log p(n)$ is positive.
   For $r=3$, it can be shown that for $n\geq 116$,
\[\Delta^3\log p(n-1)>0.\] Since
\begin{align*}
\Delta^3\log p(n-1)
=\log p(n+2)-3\log p(n+1)+3\log p(n)-\log p(n-1),
\end{align*}
it is evident that $\{p(n)\}_{n\geq116}$ is ratio log-convex. So we are led to the following assertion.

\begin{thm}\label{thsqrtpn}
The sequence $\{\sqrt[n]{p(n)}\}_{n\geq27}$ is log-convex.
\end{thm}

In the spirit of  the inequality \eqref{ineqdpc} on $\frac{p(n-1)}{p(n)}$,  we obtain the following
inequality on $\frac{\sqrt[n-1]{p(n-1)}}{\sqrt[n]{p(n)}}$.

\begin{thm}\label{th2}
For $n\geq2$, we have
\begin{equation}\label{th22}
\frac{\sqrt[n]{p(n)}}{\sqrt[n+1]{p(n+1)}}\left(1+\frac{3\pi}{\sqrt{24}n^{5/2}}\right)>\frac{\sqrt[n-1]{p(n-1)}}{\sqrt[n]{p(n)}}.
\end{equation}
\end{thm}

 Desalvo and Pak \cite{Pak} have shown that the limit of $-n^{\frac{3}{2}}\Delta^2\log p(n)$ is $\pi/\sqrt{24}$, see \eqref{limdeltapn}.
  By bounding $\Delta^2\log\sqrt[n]{p(n)}$, we derive 
  the following limit of $n^{\frac{5}{2}}\Delta^2\log\sqrt[n]{p(n)}$, which is
  analogous to  \eqref{limdeltapn}, 
\begin{equation}\label{lim3}
\lim\limits_{n \rightarrow +\infty}n^{\frac{5}{2}}\Delta^2\log\sqrt[n]{p(n)}=3\pi/\sqrt{24}.
\end{equation}
From the above relation \eqref{lim3}, it can be seen that
 the coefficent $\frac{3\pi}{\sqrt{24}}$ in \eqref{th22} is the best possible.

This paper is organized as follows. In Section 2, we show that $\{r(n)\}_{n\geq61}$ is log-convex. In Section 3, we find the limit of $n^{\frac{5}{2}}\Delta^2\log\sqrt[n]{p(n)}$ and give the  inequality \eqref{th22}.

\section{ The Log-convexity of $r(n)$}

In this section, we obtain a lower bound of $\Delta^2\log r(n-1)$ and prove the log-convexity of $\{r(n)\}_{n\geq61}$.
First, we follow the approach of Desalvo and Pak to give an expression of $\Delta^2\log r(n-1)$ as a sum of  $\Delta^2\widetilde{B}(n-1)$ and $\Delta^2\widetilde{E}(n-1)$, where $\Delta^2\widetilde{B}(n-1)$ makes a major contribution to $\Delta^2\log r(n-1)$ with $\Delta^2\widetilde{E}(n-1)$ being the error term, that is,
$\Delta^2\widetilde{B}(n-1)$ converges to $\Delta^2\log r(n-1)$.
 The expressions for $B(n)$ and $E(n)$ will be given later.
  In this setting, we derive a lower bound
 of $\Delta^2\widetilde{B}(n-1)$. By Lehmer's error bound, we
 give an upper bound for $|\Delta^2\widetilde{E}(n-1)|$.
 Combining the lower bound for $\Delta^2\widetilde{B}(n-1)$ and the upper bound for $\Delta^2\widetilde{E}(n-1)$, we are led
  to a lower bound for $\Delta^2\log r(n-1)$. By proving the positivity of this lower bound for $\Delta^2\log r(n-1)$, we reach the log-convexity of $\{r(n)\}_{n\geq61}$.

The strict log-convexity of $\{r(n)\}_{n\geq61}$ can be restated as the following relation for $n\geq61$,
\[\log r(n+1)+\log r(n-1)-2\log r(n)>0,\]
that is, for $n\geq61$,
\[\Delta^2\log r(n-1)>0.\]

For $n\geq 1$ and any positive integer $N$,
 the Hardy-Ramanujan-Rademacher formula reads
\begin{equation}\label{pnl}
p(n)=\frac{d}{\mu^2}\sum_{k=1}^{N}A_k^{\star}(n)\left[\left(1-\frac{k}{\mu}\right)e^{\frac{\mu}{k}}+\left(1+\frac{k}{\mu}\right)e^{-\frac{\mu}{k}}\right]+R_2(n,N),
\end{equation}
where $d=\frac{\pi^2}{6\sqrt{3}}$, $\mu(n)=\frac{\pi}{6}\sqrt{24n-1}$, $A_k^{\star}(n)=k^{-\frac{1}{2}}A_k(n)$, $A_k(n)$ is a   sum of 24th roots of unity with initial values $A_1(n)=1$ and $A_2(n)=(-1)^n$, $R_2(n,N)$ is the remainder.  Lehmer's error bound for $R_2(n,N)$ is given by
\begin{equation}\label{r2}
|R_2(n,N)|<\frac{\pi^2N^{-2/3}}{\sqrt{3}}\left[\left(\frac{N}{\mu}\right)^3\sinh\frac{\mu}{N}+\frac{1}{6}-\left(\frac{N}{\mu}\right)^2\right].
\end{equation}
Let us give an outline of Desalvo and Pak's approach to proving the log-concavity
of $\{p(n)\}_{n>25}$.
Setting $N=2$ in \eqref{pnl}, they expressed $p(n)$ as \begin{equation}
p(n)=T(n)+R(n),
\end{equation}
where
\begin{align}
\label{Tn} T(n)&=\frac{d}{\mu(n)^2}\left[\left(1-\frac{1}{\mu(n)}\right)e^{\mu(n)}+\frac{(-1)^n}{\sqrt{2}}e^{\frac{\mu(n)}{2}}\right],\\
\label{Rn} R(n)&=\frac{d}{\mu(n)^2}\left[\left(1+\frac{1}{\mu(n)}\right)e^{-\mu(n)}-\frac{(-1)^n}{\sqrt{2}}\frac{2}{\mu(n)}+\frac{(-1)^n}{\sqrt{2}}\left(1+\frac{2}{\mu(n)}\right)e^{-\frac{\mu(n)}{2}}\right]+R_2(n,2).
\end{align}
They  have shown that
 \begin{equation}\label{t1}
 \left|\Delta^2\log p(n-1)-\Delta^2\log T(n-1)\right|=\left|\Delta^2\log\left(1+\frac{R(n-1)}{T(n-1)}\right)\right|< e^{-\frac{\pi\sqrt{2n}}{10\sqrt{3}}}.
 \end{equation}
and
 \begin{equation}\label{t2}
 \left|\Delta^2\log T(n-1)-\Delta^2\log\frac{d}{\mu(n-1)^2}\left(1-\frac{1}{\mu(n-1)}\right)e^{\mu(n-1)}\right|<e^{-\frac{\pi\sqrt{2n}}{10\sqrt{3}}}.
 \end{equation}
 It follows
 %from \eqref{t1} and \eqref{t2}
 that $\Delta^2\log\frac{d}{\mu(n-1)^2}\left(1-\frac{1}{\mu(n-1)}\right)e^{\mu(n-1)}$ converges to $\Delta^2\log p(n-1)$. Finally, they use $-\Delta^2\log\frac{d}{\mu(n-1)^2}\left(1-\frac{1}{\mu(n-1)}\right)e^{\mu(n-1)}$ to estimate $-\Delta^2\log p(n-1)$, leading to the log-concavity of $\{p(n)\}_{n>25}$.

In this paper, we use an alternative decomposition  of $p(n)$.
 Setting $N=2$ in \eqref{pnl}, we can
   express $p(n)$ as
 \begin{equation}\label{pn}
 p(n)=\widetilde{T}(n)+\widetilde{R}(n),
 \end{equation}
 where
 \begin{equation}\label{Tbar}
\widetilde{T}(n)=\frac{d}{\mu(n)^2}\left(1-\frac{1}{\mu(n)}\right)e^{\mu(n)},
\end{equation}
\begin{align}
\label{Rnn}
\notag \widetilde{R}(n)=\frac{d}{\mu(n)^2}
&\bigg[\left(1+\frac{1}{\mu(n)}\right)
e^{-{\mu(n)}}\!+\!\frac{(-1)^n}{\sqrt{2}}\!
\left(1-\frac{2}{\mu(n)}\right)e^{\frac{\mu(n)}{2}}\\[9pt]
&+\frac{(-1)^n}{\sqrt{2}}\left(1+\frac{2}{\mu(n)}\right)e^{-\frac{\mu(n)}{2}}\bigg]
+R_2(n,2).
\end{align}
Based on the decomposition \eqref{pn} for $p(n)$,
  one can  express $\Delta^2\log r(n-1)$ as follows:
 \begin{align}\label{deltarn}
\Delta^2\log r(n-1)=\Delta^2\widetilde{B}(n-1)+\Delta^2\widetilde{E}(n-1),
 \end{align}
 where
 \begin{align}
\label{dBn}\widetilde{B}(n)&=\frac{1}{n}\log \widetilde{T}(n)-\frac{1}{n}\log n,\\[6pt]
 \label{y}\widetilde{y}_n&=\widetilde{R}(n)/\widetilde{T}(n),\\[6pt]
 \label{E}\widetilde{E}(n)&=\frac{1}{n}\log(1+\widetilde{y}_n).
\end{align}

  The following lemma will be used to give
  a lower bound and an upper bound of  $\Delta^2\widetilde{B}(n-1)$.

\begin{lem}\label{lem1}
Suppose $f(x)$ has a continuous second derivative for $x\in [n-1,n+1]$. Then there exists
$c\in (n-1,n+1)$ such that
\begin{equation}\label{i0}
\Delta^2f(n-1)=f(n+1)+f(n-1)-2f(n)=f^{''}(c).
\end{equation}
If $f(x)$ has an increasing second derivative, then
\begin{equation}\label{i1}
f''(n-1)<\Delta^2f(n-1)<f''(n+1).
\end{equation}
Conversely, if $f(x)$ has a decreasing second derivative, then
\begin{equation}\label{i2}
f''(n+1)<\Delta^2f(n-1)<f''(n-1).
\end{equation}
\end{lem}

\pf Set $\varphi(x)=f(x+1)-f(x)$. By the mean value theorem, there exists a number $\xi\in (n-1,n)$ such that
\begin{align*}
f(n+1)+f(n-1)-2f(n)=\varphi(n)-\varphi(n-1)=\varphi^{'}(\xi).
\end{align*}
Again, applying the mean value theorem to $\varphi^{'}(\xi)$, there exists a number $\theta\in(0,1)$ such that
\begin{align*}
\varphi^{'}(\xi)=f'(\xi+1)-f'(\xi)=f^{''}(\xi+\theta).
\end{align*}
Let $c=\xi+\theta$. Then we get \eqref{i0},
which yields \eqref{i1} and \eqref{i2}. \qed

In order to give a lower bound for $\Delta^2\log r(n-1)$ and obtain the limit of $n^{\frac{5}{2}}\Delta^2\log\sqrt[n]{p(n)}$,  we need the following lower and upper bounds for $\Delta^2\frac{1}{n-1}\log \widetilde{T}(n-1)$.

\begin{lem}\label{lem3}
Let
\begin{align}
B_1(n)=&\frac{72\pi}{(n+1)(24n+23)^{3/2}}
-\frac{4\log(\mu(n-1))}{(n-1)^3},\\[6pt]
B_2(n)=&\frac{72\pi}{(n-1)(24n-25)^{3/2}}
-\frac{4\log(\mu(n+1))}{(n+1)^3}+\frac{5}{(n-1)^3}.
\end{align}
For $n\geq 40$, we have
\begin{equation}
B_1(n)<\Delta^2\frac{1}{n-1}\log \widetilde{T}(n-1)<B_2(n).
\end{equation}
\end{lem}

\pf By the definition \eqref{Tbar}, we may write
\[\frac{\log\widetilde{T}(n)}{n}=\sum_{i=1}^4f_i,\]
where
\begin{align*}
f_1(n)&=\frac{\mu(n)}{n}, \\[6pt]
f_2(n)&=-\frac{3\log\mu(n)}{n}, \\[6pt] f_3(n)&=\frac{\log(\mu(n)-1)}{n},\\[6pt]
f_4(n)&=\frac{\log d}{n}.
\end{align*}
Thus
\begin{equation}\label{sumf}
\Delta^2\frac{1}{n-1}\log \widetilde{T}(n-1)=\sum_{i=1}^4\Delta^2f_i(n-1).
\end{equation}
Since
\begin{equation*}
f_1^{'''}(n)=\frac{\pi}{n(24n-1)^{3/2}}\left(-\frac{216}{n}
+\frac{864}{24n-1}+\frac{36}{n^2}-\frac{1}{n^3}\right),
\end{equation*}
we see that for $n\geq1$, $f_1^{'''}(n)<0$.
Similarly, it can be checked that for $n\geq 4$, $f_2^{'''}(n)>0$, $f_3^{'''}(n)<0$, and $f_4^{'''}(n)>0$.
Consequently, for $n\geq 4$, $f_1^{''}(n)$ and $f_3^{''}(n)$ are decreasing, whereas $f_2^{''}(n)$ and $f_4^{''}(n)$ are increasing. Using Lemma \ref{lem1}, for each $i$,  we can  get a lower bound and an upper bound for $\Delta^2f_i(n-1)$ in terms of $f_i^{''}(n-1)$ and $f_i^{''}(n+1)$. For example,
\[f_1^{''}(n+1)<\Delta^2f_1(n-1)<f_1^{''}(n-1).\]
So, by \eqref{sumf} we find that
\begin{equation}\label{lowbd}
\Delta^2\frac{1}{n-1}\log \widetilde{T}(n-1)>f_1^{''}(n+1)+f_2^{''}(n-1)+f_3^{''}(n+1)
+f_4^{''}(n-1),
\end{equation}
and
\begin{equation}
\Delta^2\frac{1}{n-1}\log \widetilde{T}(n-1)<f_1^{''}(n-1)+f_2^{''}(n+1)+f_3^{''}(n-1)
+f_4^{''}(n+1),
\end{equation}
where
\begin{align}
\label{f11}f_1^{''}(n)=&\frac{72\pi}{n(24n-1)^{3/2}}
-\frac{12\pi}{n^2(24n-1)^{3/2}}
+\frac{\pi}{3n^3(24n-1)^{3/2}},\\[6pt]
\label{f2}f_2^{''}(n)=&-\frac{6\log\mu(n)}{n^3}+\frac{72}{(24n-1)n^2}+
\frac{864}{n(24n-1)^2},\\[6pt]
\notag f_3^{''}(n)=&-\frac{4\pi^2}{(\mu(n)-1)^2(24n-1)n}
+\frac{2\log(\mu(n)-1)}{n^3}\\[6pt]
\label{f3}&-\frac{4\pi}{(\mu(n)-1)\sqrt{24n-1}n^2}
-\frac{24\pi}{(\mu(n)-1)(24n-1)^{3/2}n},\\[6pt]
\label{f4}f_4^{''}(n)=&\frac{2\log d}{n^3}.
\end{align}
According to \eqref{f11}, one can check that  for $n\geq2$,
\begin{equation}\label{f1}
f_1^{''}(n+1)>\frac{72\pi}{(n+1)(24n+23)^{3/2}}
-\frac{12\pi}{(n+1)^2(24n+23)^{3/2}}.
\end{equation}
An easy computation shows that for $n\geq3$,
\begin{equation}\label{23mu}
\mu(n)-1>\frac{2}{3}\mu(n-2).
\end{equation}
 Substituting  \eqref{23mu} into \eqref{f3} yields that
\begin{equation}\label{f3b}
f_3^{''}(n+1)>\frac{2\log(\mu(n+1)-1)}{(n+1)^3}
-\frac{540}{(24n-25)^2(n-1)}-\frac{36}{(24n-25)(n-1)^2}.
\end{equation}
Using \eqref{f2} and \eqref{f3b}, we find that
\begin{align}\label{f2f3}
\notag &f_2^{''}(n-1)+f_3^{''}(n+1)\\[6pt]
\notag&\qquad>\frac{2\log(\mu(n+1)-1)}{(n+1)^3}-\frac{6\log(\mu(n-1))}{(n-1)^3}\\[6pt]
&\qquad~~~~~~ +\frac{324}{(n-1)(24n-25)^2}+\frac{36}{(n-1)^2(24n-25)}
\end{align}
Apparently,   for $n\geq 2$,
\begin{align*}
\frac{2}{(n+1)^3}-\frac{2}{(n-1)^3}>-\frac{12}{(n-1)^4},
\end{align*}
so that
\begin{align}\label{f26}
\notag &\frac{2\log(\mu(n+1)-1)}{(n+1)^3}
-\frac{6\log(\mu(n-1))}{(n-1)^3}\\[6pt]
\notag&\qquad >\frac{2\log(\mu(n+1)-1)}{(n+1)^3}-\frac{2\log(\mu(n+1)-1)}{(n-1)^3}
-\frac{4\log(\mu(n-1))}{(n-1)^3}\\[6pt]
&\qquad >-\frac{12\log(\mu(n+1)-1)}{(n-1)^4}-\frac{4\log(\mu(n-1))}{(n-1)^3}.
\end{align}
Since, for $n\geq2$,  
\begin{equation}\label{f231}
\frac{324}{(n-1)(24n-25)^2}+\frac{36}{(n-1)^2(24n-25)}
>\frac{2}{(n-1)^3},
\end{equation}
utilizing  \eqref{f2f3} and \eqref{f26} yields that for $n\geq 3$,
\begin{align}\label{f232}
f_2^{''}(n-1)+f_3^{''}(n+1)>-\frac{4\log(\mu(n-1))}{(n-1)^3}+\frac{
2}{(n-1)^3}-\frac{12\log(\mu(n+1)-1)}{(n-1)^4}.
\end{align}
Using \eqref{f4}, \eqref{f1} and \eqref{f232}, we deduce that
\begin{align}
\notag &f_1^{''}(n+1)+f_2^{''}(n-1)+f_3^{''}(n+1)+f_4^{''}(n-1)-B_1(n)\\[6pt]
\label{Cn}&\qquad>\frac{2(1+\log d)}{(n-1)^3}-\frac{12\pi}{(n+1)^2(24n+23)^{3/2}}
-\frac{12\log(\mu(n+1)-1)}{(n-1)^4}.
\end{align}
Let $C(n)$ be the right hand side of \eqref{Cn}. 
To prove \eqref{lowbd}, it is enough to show that $C(n)>0$ when $n\geq40$.
Since $\log x<x$ for $x>0$, and for $n\geq3$
\begin{align}\label{mu}
\mu(n+1)-1<\frac{\pi}{4}\sqrt{24n-24},
\end{align}
 we get
\begin{equation}\label{cn1}
-\frac{12\log(\mu(n+1)-1)}{(n-1)^4}>-\frac{12(\mu(n+1)-1)}{(n-1)^4}>-\frac{
3\sqrt{24}\pi}{(n-1)^{7/2}}.
\end{equation}
Note that for $n\geq2$,
\begin{equation}\label{cn2}
-\frac{12\pi}{(n+1)^2(24n+23)^{3/2}}>
-\frac{\sqrt{24}\pi}{48(n-1)^{7/2}}.
\end{equation}
Combining \eqref{cn1} and \eqref{cn2} gives for $n\geq 2$,
\begin{equation}\label{d}
C(n)>\frac{2(1+\log d)}{(n-1)^3}-\frac{(3+1/48)\sqrt{24}\pi}{(n-1)^{7/2}}.
\end{equation}
It is straightforward to show that the right
hand side of \eqref{d} is positive if $n\geq 490$.
For $40\leq n\leq 489$, it is routine to check that
$C(n)>0$, and so $C(n)>0$ for $n\geq 40$.
It follows from \eqref{Cn} that for  $n\geq40$,
\[\Delta^2\frac{1}{n-1}\log \widetilde{T}(n-1)>B_1(n).\]
To derive the upper bound for $\Delta^2\frac{1}{n-1}\log \widetilde{T}(n-1)$, we obtain the following upper bounds
 which can be verified directly. The proofs are omitted.
  For $n\geq2$,  
\begin{align*}
f_1^{''}(n-1)<&\frac{72\pi}{(n-1)[24n-25]^{3/2}},\\[6pt]
f_2^{''}(n+1)<&-\frac{6\log\mu(n+1)}{(n+1)^3}+\frac{9}{2(n-1)^3},\\[6pt]
f_3^{''}(n-1)<&-\frac{4\pi^2}{(\mu(n-1))^2(24n-25)(n-1)}
+\frac{2\log(\mu(n-1))}{(n-1)^3}\\[6pt]
&\qquad -\frac{4\pi}{\mu(n-1)\sqrt{24n-25}(n-1)^2}
-\frac{24\pi}{\mu(n-1)(24n-25)^{3/2}(n-1)},\\[6pt]
f_2^{''}(n+1)&+f_3^{''}(n-1)<\frac{3}{(n-1)^3}
+\frac{12\log(\mu(n+1))}{(n-1)^4}-\frac{4\log(\mu(n+1))}{(n+1)^3},\\[6pt]
&\qquad \qquad \qquad \qquad f_4^{''}(n+1)<0.
\end{align*}
Combining the above upper bounds, we conclude that for $n\geq40$,
\begin{equation*}
f_1^{''}(n-1)+f_2^{''}(n+1)+f_3^{''}(n-1)
+f_4^{''}(n+1)<B_2(n).
\end{equation*}
This completes the proof. \qed\\

The following lemma gives an upper bound for $|\Delta^2 \widetilde{E}(n-1)|$.

\begin{lem}\label{lem2}
For $n\geq 40$,
\begin{equation}\label{lem2ie}
|\Delta^2\widetilde{E}(n-1)|<
\frac{5}{n-1}e^{-\frac{\pi\sqrt{24n-25}}{18}}.
\end{equation}
\end{lem}

\pf By \eqref{E}, we find that for $n\geq 2$,
\begin{align}\label{deltaen}
\Delta^2\widetilde{E}(n-1)=\frac{1}{n-1}
\log(1+\widetilde{y}_{n-1})+\frac{1}{n+1}\log(1+\widetilde{y}_{n+1})
-\frac{2}{n}\log(1+\widetilde{y}_n),
\end{align}
where
\[\widetilde{y}_n=\widetilde{R}(n)/\widetilde{T}(n).\]
To bound $|\Delta^2\widetilde{E}(n-1)|$, it is necessary to bound $\widetilde{y}_n$. For this purpose, we first consider $\widetilde{R}(n)$,  as defined by \eqref{Rnn}. Since $d<1$ and $\mu(n)>2$,   for $n\geq1$ we have
\begin{align*}
\frac{d}{\mu(n)^2}&\left[\left(1\!+\!\frac{1}{\mu(n)}\right)\!
e^{-{\mu(n)}}\!+\!\frac{(-1)^n}{\sqrt{2}}\!
\left(1\!-\!\frac{2}{\mu(n)}\right)\!e^{\frac{\mu(n)}{2}}\!
+\!\frac{(-1)^n}{\sqrt{2}}\!\left(1\!+\!\frac{2}{\mu(n)}\right)\!
e^{-\frac{\mu(n)}{2}}\right]\\[6pt]
&<\frac{1}{\mu(n)^2}\left(1\!+\!e^{\frac{\mu(n)}{2}}\!+\!1\right).
\end{align*}
For $N=2$ and $n\geq 1$, Lehmer's bound \eqref{r2} reduces to
\begin{equation*}
|R_2(n,2)|<4\left(1+\frac{4}{\mu(n)^3}e^{\frac{\mu(n)}{2}}\right).
\end{equation*}
By the definition of $\widetilde{R}(n)$,   
\begin{align}\label{Rnleq}
|\widetilde{R}(n)|<\frac{1}{\mu(n)^2}\left(1+e^{\frac{\mu(n)}{2}}
+1\right)+4\left(1+\frac{4}{\mu(n)^3}e^{\frac{\mu(n)}{2}}\right)
<5+\frac{9}{\mu(n)^2}e^{\frac{\mu(n)}{2}}.
\end{align}
Recalling the definition \eqref{Tbar} of $\widetilde{T}(n)$,
it follows from \eqref{Rnleq} that for $n\geq1$,
\begin{equation}\label{yn}
|\widetilde{y}_n|<\frac{\mu(n)}{d(\mu(n)-1)}\left(5\mu(n)^2
e^{-\frac{2\mu(n)}{3}}+9e^{-\frac{\mu(n)}{6}}\right)e^{-\frac{\mu(n)}{3}}.
\end{equation}
Observe that for $n\geq2$,
\begin{align}\label{de1}
\left(5\mu(n)^2e^{-\frac{2\mu(n)}{3}}+9e^{-\frac{\mu(n)}{6}}\right)^{'}<0,
\end{align}
 and
\begin{align}\label{de2}
 \left( \frac{d(\mu(n)-1)}{\mu(n)}\right)^{'}>0.
\end{align}
Since
\begin{align*}
5\mu^2(40)e^{-\frac{2\mu(40)}{3}}
+9e^{-\frac{\mu(40)}{6}}<\frac{d(\mu(40)-1)}{\mu(40)},
\end{align*}
  using \eqref{de1} and \eqref{de2},  we deduce that for $n\geq40$,
\begin{equation}\label{yn1}
5\mu^2(n)e^{-\frac{2\mu(n)}{3}}+9e^{-\frac{\mu(n)}{6}}<
\frac{d(\mu(n)-1)}{\mu(n)}.
\end{equation}
Now, it is clear from (\ref{yn}) and (\ref{yn1})  that for $n\geq 40$,
\begin{equation}\label{yn2}
|\widetilde{y}_n|<e^{-\frac{\mu(n)}{3}}.
\end{equation}
In view of \eqref{yn2},  for $n\geq40$,
\begin{equation}\label{yn4}
|\widetilde{y}_n|<e^{-\frac{\mu(40)}{3}}<\frac{1}{5}.
\end{equation}
It is known that  $\log(1+x)<x$ for $0<x<1$ and $-\log(1+x)<-x/(1+x)$ for $-1<x<0$. Thus, for $|x|<1$,
\begin{equation}\label{inlog}
|\log(1+x)|\leq\frac{|x|}{1-|x|},
\end{equation}
see also \cite{Pak}, and so it follows from
  \eqref{yn4} and \eqref{inlog}  that for $n\geq40$,
\begin{equation}\label{yn3}
|\log(1+\widetilde{y}_n)|\leq
\frac{|\widetilde{y}_n|}{1-|\widetilde{y}_n|}\leq
\frac{5}{4}|\widetilde{y}_n|.
\end{equation}
Because of \eqref{deltaen}, we see that for $n\geq2$,
\begin{align}\label{deltaen1}
\left|\Delta^2\widetilde{E}(n\!-\!1)\right|
\leq\frac{1}{n\!-\!1}\left|\log(1\!+\!\widetilde{y}_{n-1})\right|
\!+\!\frac{1}{n\!+\!1}\left|\log(1\!+\!\widetilde{y}_{n+1})\right|\!+
\!\frac{2}{n}\left|\log(1\!+\!\widetilde{y}_n)\right|.
\end{align}
Applying \eqref{yn3} to \eqref{deltaen1}, we obtain that for
  $n\geq40$,
\begin{align}\label{delataen2}
\left|\Delta^2\widetilde{E}(n-1)\right|
\leq\frac{5}{4}\left[\frac{|\widetilde{y}_{n-1}|}{n-1}+
\frac{|\widetilde{y}_{n+1}|}{n+1}+\frac{2|\widetilde{y}_n|}{n}\right].
\end{align}
Plugging (\ref{yn2}) into \eqref{delataen2}, we infer that for $n\geq40$,
\begin{align}\label{in8}
\left|\Delta^2\widetilde{E}(n-1)\right|<\frac{5}{4}
\left[\frac{e^{-\frac{\mu(n-1)}{3}}}{n-1}+
\frac{e^{-\frac{\mu(n+1)}{3}}}{n+1}+
\frac{2e^{-\frac{\mu(n)}{3}}}{n}\right].
\end{align}
But
$\frac{1}{n}e^{-\frac{\mu(n)}{3}}$ is decreasing for $n\geq 1$, it follows from \eqref{in8} that for $n\geq40$,
\[\left|\Delta^2\widetilde{E}(n-1)\right|<
\frac{5}{n-1}e^{-\frac{\mu(n-1)}{3}}.\]
This proves \eqref{lem2ie}.\qed\\

With the aid of Lemmas Lemma \ref{lem3} and \ref{lem2}, we are ready to prove the log-convexity of $\{r(n)\}_{n\geq61}$.

\noindent{\it Proof of Theorem \ref{th1}}. To prove the strict log-convexity of $\{r(n)\}_{n\geq61}$, we proceed to show that  for $n\geq61$,
\[\Delta^2\log r(n-1)>0.\]
Evidently, for $n\geq40$,
\[\left(-\frac{\log n}{n}\right)'''>0. \]
By Lemma \ref{lem1},
\[-\Delta^2\frac{\log(n-1)}{n-1}>\left(-\frac{\log (n-1)}{n-1}\right)'',\]
that is,
\begin{equation}\label{logn}
-\Delta^2\frac{\log(n-1)}{n-1}>-\frac{2\log(n-1)}{(n-1)^3}
+\frac{1}{(n-1)^3}.
\end{equation}
It follows from \eqref{dBn} that
\[\Delta^2\widetilde{B}(n-1)=\Delta^2\frac{1}{n-1}\log \widetilde{T}(n-1)-\Delta^2\frac{\log(n-1)}{n-1}.\]
Applying Lemma \ref{lem3} and \eqref{logn} to the above relation, we deduce that for $n\geq40$,
\[\Delta^2\widetilde{B}(n-1)>\widetilde{B}_1(n)-
\frac{2\log(n-1)}{(n-1)^3}+\frac{3}{(n-1)^3},\]
that is,
\begin{equation}\label{db1}
\Delta^2\widetilde{B}(n-1)>\frac{72\pi}{(n+1)[24n+23]^{3/2}}-
\frac{4\log[\mu(n-1)]}{(n-1)^3}
-\frac{2\log(n-1)}{(n-1)^3}+\frac{3}{(n-1)^3}.
\end{equation}
By \eqref{deltarn} and Lemma \ref{lem2}, we find that for $n\geq40$,
\begin{align}\label{deltarn1}
\Delta^2\log r(n-1)
>\Delta^2\widetilde{B}(n-1)-\frac{5}{n-1}
e^{-\frac{\pi\sqrt{24n-25}}{18}}.
\end{align}
It follows from \eqref{db1} and \eqref{deltarn1}  that for $n\geq40$,
\begin{align*}
&\Delta^2\log r(n-1)\\[6pt]
&>\frac{72\pi}{(n+1)[24n+23]^{3/2}}-\frac{4\log[\mu(n-1)]}{(n-1)^3}
-\frac{2\log(n-1)}{(n-1)^3}+\frac{3}{(n-1)^3}-\frac{5}{n-1}
e^{-\frac{\pi\sqrt{24n-25}}{18}}.
\end{align*}
Let $D(n)$ denote the right hand side of the above relation.
Clearly, for $n\geq5505$,
\begin{align}\label{c4}
\frac{72\pi}{(n+1)[24n+23]^{3/2}}>\frac{3\pi}{\sqrt{24}(n+1)^{5/2}}
>\frac{1}{(n-1)^{5/2}}.
\end{align}
To prove that $D(n)>0$ for $n\geq 5505$, we wish to show that for $n\geq5505$,
\begin{equation}\label{c}
-\frac{4\log[\mu(n-1)]}{(n-1)^3}
-\frac{2\log(n-1)}{(n-1)^3}+\frac{3}{(n-1)^3}-\frac{5}{n-1}
e^{-\frac{\pi\sqrt{24n-25}}{18}}>-\frac{1}{(n-1)^{5/2}}.
\end{equation}
Using the fact that for $x>5504$, $\log x<x^{1/4}$,
we deduce that for $n\geq5505$,
\begin{align}\label{c1}
\frac{4\log[\mu(n-1)]}{(n-1)^3}<\frac{4\sqrt[4]{\mu(n-1)}}{(n-1)^3}<
\frac{4\sqrt[4]{\frac{\pi}{4}\sqrt{24n-24}}}{(n-1)^3}
<\frac{6}{(n-1)^{23/8}},
\end{align}
and
\begin{align}\label{c2}
\frac{2\log(n-1)}{(n-1)^3}<\frac{2(n-1)^{1/4}}{(n-1)^3}
<\frac{2}{(n-1)^{11/4}}.
\end{align}
Since $e^x>x^6/720$ for $x>0$,  we see that for $n\geq2$,
\begin{align}\label{ex}
\frac{1}{n-1}e^{-\frac{\pi\sqrt{24n-25}}{18}}
<\frac{1}{n-1}e^{-\frac{\pi\sqrt{23n}}{18}}
<\frac{2094}{n^3(n-1)}<\frac{2094}{(n-1)^4}.
\end{align}
Combining \eqref{c1}, \eqref{c2} and \eqref{ex}, we find that for $n\geq5505$,
 \begin{align*}
 &-\frac{4\log[\mu(n-1)]}{(n-1)^3}
-\frac{2\log(n-1)}{(n-1)^3}+\frac{3}{(n-1)^3}-\frac{5}{n-1}
e^{-\frac{\pi\sqrt{24n-25}}{18}}\\
 &\qquad >-\frac{6}{(n-1)^{23/8}}-\frac{2}{(n-1)^{11/4}}+\frac{3}{(n-1)^3}
 -\frac{10470}{(n-1)^4}\\
 &\qquad >-\frac{6}{(n-1)^{23/8}}-\frac{2}{(n-1)^{11/4}}\\
 &\qquad>-\frac{1}{(n-1)^{5/2}}.
 \end{align*}
This proves the  inequality \eqref{c}. By \eqref{c} and \eqref{c4}, we obtain that $D(n)>0$ for $n\geq5505$.
Verifying  that
 $\Delta^2\log r(n-1)$ for $61\leq n\leq5504$ completes the proof. \qed

\section{An inequality on the ratio $\frac{\sqrt[n-1]{p(n-1)}}{\sqrt[n]{p(n)}}$ }

In this section, we employ Lemma \ref{lem3} and Lemma \ref{lem2} to find the limit of $n^{\frac{5}{2}}\Delta^2\log\sqrt[n]{p(n)}$. Then we give an upper bound for $\Delta^2\log\sqrt[n-1]{p(n-1)}$. This leads to an inequality analogous to the inequality \eqref{ineqdpc}.

\begin{thm}\label{lem4}
Let $\alpha=3\pi/\sqrt{24}$. We have
\begin{equation}
\lim\limits_{n \rightarrow +\infty}n^{\frac{5}{2}}\Delta^2\log\sqrt[n]{p(n)}=\alpha.
\end{equation}
\end{thm}

\pf Using \eqref{pn}, that is, the $N=2$ case of the Hardy-Ramanujan-Rademacher formula  for $p(n)$,
we find that
\begin{align*}
\log\sqrt[n]{p(n)}&=\frac{1}{n}\log[\widetilde{T}(n)+\widetilde{R}(n)]\\
&=\frac{1}{n}\log \widetilde{T}(n)+\frac{1}{n}\log\left(1+\frac{\widetilde{R}(n)}{\widetilde{T}(n)}\right)\\[6pt]
&=\frac{1}{n}\log \widetilde{T}(n)+\frac{1}{n}\log(1+\widetilde{y}_n),
\end{align*}
where $\widetilde{T}(n)$ and $y_n$ are given by \eqref{Tbar} and \eqref{y}.
By the definition \eqref{E} of $\widetilde{E}(n)$, we get
\begin{equation}\label{lim0}
\Delta^2\log\sqrt[n-1]{p(n-1)}=\Delta^2\frac{1}{n-1}\log \widetilde{T}(n-1)+\Delta^2\widetilde{E}(n-1).
\end{equation}
Applying Lemma \ref{lem3}, we obtain that  for $n\geq 40$,
\begin{equation}\label{b}
B_1(n)<\Delta^2\frac{1}{n-1}\log \widetilde{T}(n-1)<B_2(n),
\end{equation}
where
\begin{align*}
B_1(n)=&\frac{72\pi}{(n+1)[24n+23]^{3/2}}-
\frac{4\log[\mu(n-1)]}{(n-1)^3},\\[6pt]
B_2(n)=&\frac{72\pi}{(n-1)[24n-25]^{3/2}}-
\frac{4\log[\mu(n+1)]}{(n+1)^3}+\frac{5}{(n-1)^3}.
\end{align*}
It is easily seen that
\begin{align}
\label{limb1}&\lim\limits_{n \rightarrow +\infty}\frac{72\pi(n-1)^{5/2}}{(n+1)[24n+23]^{3/2}}=\alpha,\\[6pt] \label{limb2}&\lim\limits_{n \rightarrow +\infty}\frac{\log\mu}{(n-1)^{1/2}}=0.
\end{align}
By \eqref{limb1} and \eqref{limb2}, we see that
\begin{equation}\label{b1}
\lim\limits_{n\rightarrow+\infty}(n-1)^{\frac{5}{2}}B_1(n)
=\lim\limits_{n\rightarrow+\infty}(n-1)^{\frac{5}{2}}B_2(n)=\alpha.
\end{equation}
Combining \eqref{b} and \eqref{b1} gives
\begin{align}\label{limt}
\lim\limits_{n\rightarrow+\infty}(n-1)^{\frac{5}{2}}\Delta^2\frac{1}{n-1}\log \widetilde{T}(n-1)=\alpha.
\end{align}
 From Lemma \ref{lem2}, we know that for $n\geq 40$,
\begin{equation*}
-\frac{5}{n-1}e^{-\frac{\pi\sqrt{24n-25}}{18}}
<\Delta^2\widetilde{E}(n-1)<\frac{5}{n-1}
e^{-\frac{\pi\sqrt{24n-25}}{18}}.
\end{equation*}
By the fact that
 \[\lim\limits_{n \rightarrow +\infty}(n-1)^{\frac{3}{2}}e^{-\frac{\pi\sqrt{24n}}{18}} =0,\]
we get
\begin{align}\label{lime}
\lim\limits_{n\rightarrow+\infty}(n-1)^{\frac{5}{2}}\Delta^2\widetilde{E}(n-1)=0.
\end{align}
Using \eqref{lim0}, \eqref{limt} and \eqref{lime}, we deduce that
\begin{equation*}
\lim\limits_{n \rightarrow +\infty}n^{\frac{5}{2}}\Delta^2\log\sqrt[n]{p(n)}=\alpha,
\end{equation*}
as required. \qed

To prove Theorem \ref{th2}, we
 need the following   upper bound for $\Delta^2\log\sqrt[n-1]{p(n-1)}$.

\begin{thm}\label{lem5}
For $n\geq 2095$,
\begin{equation}\label{3.9}
\Delta^2\log\sqrt[n-1]{p(n-1)}<\frac{3\pi}{\sqrt{24}n^{5/2}+3\pi}.
\end{equation}
\end{thm}
\pf By the upper bound of $\Delta^2\frac{1}{n-1}\log \widetilde{T}(n-1)$  given in Lemma \ref{lem3},
the upper bound of $\Delta^2\widetilde{E}(n-1)$ given in Lemma \ref{lem2} and the relation \eqref{lim0},
we get the following upper bound of $\Delta^2\log\sqrt[n-1]{p(n-1)}$ for $n\geq40$,
 \begin{equation*}
 \Delta^2\log\sqrt[n-1]{p(n-1)}<\frac{72\pi}{(n-1)[24n-25]^{3/2}}+\frac{5}{(n-1)^3}-
\frac{4\log[\mu(n+1)]}{(n+1)^3}+\frac{5}{n-1}
e^{-\frac{\pi\sqrt{24n-25}}{18}}.
 \end{equation*}
To prove \eqref{3.9}, we claim that for $n\geq 2095$,
\begin{equation}\label{cla}
\frac{72\pi}{(n-1)[24n-25]^{3/2}}+\frac{5}{(n-1)^3}-
\frac{4\log[\mu(n+1)]}{(n+1)^3}+\frac{5}{n-1}
e^{-\frac{\pi\sqrt{24n-25}}{18}}<\frac{3\pi}{\sqrt{24}n^{5/2}+3\pi}.
\end{equation}
First, we show that for $n\geq60$,
\begin{equation}\label{g1}
\frac{72\pi}{(n-1)[24n-25]^{3/2}}-\frac{3\pi}{\sqrt{24}n^{5/2}+3\pi}
<\frac{1}{(n-1)^3}.
\end{equation}
For $0<x\leq\frac{1}{48}$, it can be checked that
\begin{equation}\label{x}
\frac{1}{(1-x)^{3/2}}<1+\frac{3}{2}x+\frac{3}{8}x^{\frac{3}{2}}.
\end{equation}
In the notation $\alpha=3\pi/\sqrt{24}$, we have
\begin{align}\label{xx}
\frac{72\pi}{(n-1)(24n-25)^{3/2}}=\frac{\alpha}{(n-1)n^{3/2}
(1-\frac{25}{24n})^{3/2}}.
\end{align}
 Setting $x=\frac{25}{24n}$, we have  $x \leq\frac{1}{48}$ for $n\geq60$. Applying \eqref{x} to the right hand side of \eqref{xx}, we find that for $n
 \geq 60$,
\begin{align}\label{g3}
\frac{\alpha}{(n-1)n^{3/2}(1-\frac{25}{24n})^{3/2}}
<\frac{\alpha}{(n-1)n^{3/2}}\left[1+\frac{75}{48n}
+\frac{3}{8}\left(\frac{25}{24n}\right)^{\frac{3}{2}}\right],
\end{align}
so  that for $n\geq 60$,
\begin{align}
\notag &\frac{72\pi}{(n-1)[24n-25]^{3/2}}-\frac{3\pi}{\sqrt{24}n^{5/2}+3\pi}\\[6pt] \label{g5}&\qquad<\frac{\alpha}{(n-1)n^{3/2}}-\frac{3\pi}{\sqrt{24}n^{5/2}+3\pi}
+\frac{\alpha}{(n-1)n^{3/2}}\left[\frac{75}{48n}
+\frac{3}{8}\left(\frac{25}{24n}\right)^{\frac{3}{2}}\right].
\end{align}
To prove \eqref{g1}, we proceed to show that the
right hand side of \eqref{g5}   is bounded by
${1\over (n-1)^3}$.
For $n\geq2$, we obtain that
\begin{align*} &\frac{\alpha}{(n-1)n^{3/2}}-\frac{3\pi}{\sqrt{24}n^{5/2}+3\pi}\\[6pt]
&\qquad=\frac{\alpha}{(n-1)n^{3/2}}-\frac{\alpha}{n^{5/2}+\alpha}\\[6pt]
&\qquad=\frac{\alpha n^{3/2}+\alpha^2}{(n^{5/2}+\alpha)(n-1)n^{3/2}}\\[6pt]
&\qquad=\frac{\alpha}{(n^{5/2}+\alpha)(n-1)}+\frac{\alpha^2}{(n^{5/2}
+\alpha)(n-1)n^{3/2}}.
\end{align*}
Since $n^{5/2}+\alpha>(n-1)^{5/2}$ and $n^{3/2}>(n-1)^{3/2}$ for $n\geq 2$, in this case we have
\begin{align}\label{g4}
\frac{\alpha}{(n-1)n^{3/2}}-\frac{3\pi}{\sqrt{24}n^{5/2}+3\pi}
<\frac{\alpha}{(n-1)^{7/2}}+\frac{\alpha}{(n-1)^5}.
\end{align}
Applying \eqref{g4} to \eqref{g5}, we obtain that for $n\geq60$,
\begin{align}
\notag &\frac{72\pi}{(n-1)[24n-25]^{3/2}}-\frac{3\pi}{\sqrt{24}n^{5/2}+3\pi}\\[6pt] \label{g6}&\qquad<\frac{\alpha}{(n-1)^{7/2}}+\frac{\alpha^2}{(n-1)^5}
+\frac{\alpha}{(n-1)n^{3/2}}\left[\frac{75}{48n}
+\frac{3}{8}\left(\frac{25}{24n}\right)^{\frac{3}{2}}\right].
\end{align}
Since $\frac{75}{48n}<\frac{2}{n-1}$ and $\frac{3}{8}\left(\frac{25}{24n}\right)^{\frac{3}{2}}
<\frac{1}{(n-1)^{3/2}}$ for $n\geq 2$,
 it follows from \eqref{g6}  that for $n\geq 60$,
\begin{align*}
&\frac{72\pi}{(n-1)[24n-25]^{3/2}}-\frac{3\pi}{\sqrt{24}n^{5/2}+3\pi}\\[6pt] &\qquad<\frac{\alpha}{(n-1)^{7/2}}+\frac{\alpha^2}{(n-1)^5}
+\frac{2\alpha}{(n-1)^{7/2}}+\frac{\alpha}{(n-1)^4}.
\end{align*}
Using the fact that $\alpha<2$, we see that
\begin{align}\label{g7}
\frac{3\alpha}{(n-1)^{7/2}}+\frac{\alpha^2}{(n-1)^5}
+\frac{\alpha}{(n-1)^4}
< \frac{6}{(n-1)^{7/2}}+\frac{4}{(n-1)^5}+\frac{2}{(n-1)^4}.
\end{align}
For $n\geq 60$, it is easily checked that the right hand side of \eqref{g7} is bounded by $\frac{1}{(n-1)^3}$.
This confirms \eqref{g1}.

 To prove the claim \eqref{cla},
 it is enough to show that for $n\geq2095$,
\begin{equation}\label{g2}
\frac{1}{(n-1)^3}<\frac{4\log[\mu(n+1)]}{(n+1)^3}-\frac{5}{(n-1)^3}
-\frac{5}{n-1}e^{-\frac{\pi\sqrt{24n-25}}{18}}.
\end{equation}
From  \eqref{ex}  it can be seen that for $n\geq2095$,
\begin{align}\label{g8}
\frac{5}{n-1}e^{-\frac{\pi\sqrt{24n-25}}{18}}<\frac{5}{(n-1)^3}.
\end{align}
Since $4\log[\mu(n+1)]>18$ for $n\geq2095$,
it follows from \eqref{g8} that for $n\geq2095$,
\begin{align*}
&\frac{4\log[\mu(n+1)]}{(n+1)^3}-\frac{5}{(n-1)^3}
-\frac{5}{n-1}e^{-\frac{\pi\sqrt{24n-25}}{18}}\\[6pt]
&\qquad >\frac{18}{(n+1)^3}-\frac{10}{(n-1)^3}>\frac{1}{(n-1)^3}.
\end{align*}
So we obtain \eqref{g2}. 
  This completes the proof. \qed

We are now in a position to finish the proof of Theorem \ref{th2}.

\noindent{\it Proof of Theorem \ref{th2}}.  It is known that for $x>0$,
\[\frac{x}{1+x}<\log(1+x),\]
so  that for $n\geq 1$,
\[\frac{3\pi}{\sqrt{24}n^{5/2}+3\pi}
<\log\left(1+\frac{3\pi}{\sqrt{24}n^{5/2}}\right).\]
In light of  the above relation, Theorem \ref{lem5} implies that for $n\geq2095$,
\begin{equation*}
\Delta^2\log\sqrt[n-1]{p(n-1)}
<\log\left(1+\frac{3\pi}{\sqrt{24}n^{5/2}}\right),
\end{equation*}
that is,
 \[\sqrt[n+1]{p(n+1)}\sqrt[n-1]{p(n-1)}
 <\left(1+\frac{3\pi}{\sqrt{24}n^{5/2}}\right)(\sqrt[n]{p(n)})^2.\] It can be checked that  the above inequality holds for $2\leq n\leq2095$.  This completes the proof of the theorem.\qed

We remark that $\alpha=3\pi/\sqrt{24}$ is the smallest possible number for the inequality in Theorem \ref{th2}. Suppose that $0< \beta<\alpha$.
By Theorem \ref{lem4}, there exists an integer $N$ so as to for $n>N$,
\[n^{5/2}\Delta^2\log \sqrt[n-1]{p(n-1)}>\beta.\] It follows that
\[\Delta^2\log \sqrt[n-1]{p(n-1)}>\frac{\beta}{n^{5/2}}
>\log\left(1+\frac{\beta}{n^{5/2}}\right),\]
which implies that  for $n>N$,
\[\frac{\sqrt[n]{p(n)}}{\sqrt[n+1]{p(n+1)}}
\left(1+\frac{\beta}{n^{5/2}}\right)
<\frac{\sqrt[n-1]{p(n-1)}}{\sqrt[n]{p(n)}}.\]

\vspace{0.5cm}
 \noindent{\bf Acknowledgments.} This work was supported by  the 973
Project, the PCSIRT Project of the Ministry of Education,  and the National Science
Foundation of China.


\begin{thebibliography}{99}

\bibitem{Andrews1}
G.E. Andrews, Combinatorial proof of a partition function limit, Amer. Math. Monthly., 78 (1971), 276-278.

%\bibitem{Andrews2}
%G.E. Andrews, The Theory of Partitions, Cambridge University Press, Cambridge (1998)

\bibitem{Chen}
W.Y.C. Chen, Recent developments on log-concavity and q-log-concavity of combinatorial polynomials, DMTCS Proceedings of 22nd International Conference on Formal Power Series and Algebraic Combinatorics, 2010.

\bibitem{Chen-a}
W.Y.C. Chen, J.J.F. Guo and L.X.W. Wang, Infinitely-logarithmically monotonic combinatorial sequences. Adv. Appl. Math., 52 (2014), 99-120.

\bibitem{Chen1}
W.Y.C. Chen, L.X.W. Wang and G.Y.B. Xie, Finite differences of the logarithm of the partition function. Math. Comp., to appear.

\bibitem{Pak}
S. Desalvo and I. Pak, Log-concavity of the partition function,
Ramanujan J., to appear.

\bibitem{Hardy1}
G.H. Hardy, Twelve Lectures on Subjects Suggested by His Life and Work, Cambridge University Press,
Cambridge, 1940.

\bibitem{Hardy}
G.H. Hardy and S. Ramanujan, Asymptotic formulae in combinatory
analysis, Proc. London Math. Soc., 17 (1918), 75-175.

\bibitem{Leh1}
D.H. Lehmer, On the series for the partition function, Trans. Amer. Math. Soc., 43 (1938), 271-292.

\bibitem{Leh2}
 D.H. Lehmer, On the remainders and convergence of the series for the partition function, Trans. Amer. Math. Soc., 46 (1939), 362-373.

\bibitem{Rade}
H. Rademacher, A convergent series for the partition function $p(n)$,
Proc. Nat. Acad. Sci, 23 (1937), 78-84.

\bibitem{Sun}
 Z.-W. Sun, On a sequence involving sums of primes, Bull. Aust. Math. Soc. 88 (2013), 197-205.

\bibitem{Wang}
Y. Wang, B.-X. Zhu, Proofs of some conjectures on monotonicity of number-theoretic and combinatorial sequences, Sci. China Math. 57 (2014), 2429-2435.







































\end{thebibliography}
\end{document}